\documentclass{article}

\usepackage{algpseudocode}
\usepackage{amsfonts}
\usepackage{amsmath}
\usepackage{booktabs}

\newtheorem{algorithm}{Algorithm}
\newtheorem{lemma}{Lemma}
\def\Re{\mathbb{R}}
\def\bu{\mathbf{u}}
\def\ulp{\mbox{\tt ulp}}
\def\Nan{\mbox{\tt NaN}}
\def\d1{\delta_1}

\title{A Square-Root Free Algorithm for Computing Real Givens Rotations}

\author{Carlos F. Borges}

\begin{document}

\maketitle

% \renewcommand{\thefootnote}{\fnsymbol{footnote}}

% \footnotetext[2]{Department of Applied Mathematics, Naval Postgraduate School, Monterey CA 93943. borges@nps.edu}

\begin{abstract}
   We develop an accurate square-root-free algorithm for constructing real Givens rotations. On processors that support the fused multiply-add operation in hardware, the algorithm is competitive with square-root based algorithms using a hardware square-root. Unlike the square-root-free algorithms in \cite{Hsieh1993,GENTLEMAN1973,Barlow1987,Hammarling1974ANO,Ling1989EfficientLL,Hanson,Hanson2}, our approach will construct the Givens rotation directly and is therefore applicable to a much wider variety of algorithms that use Givens rotations. We investigate the accuracy of the algorithm by simulation.
\end{abstract}

% \begin{keywords}
% Givens rotation, floating-point, fused multiply-add, IEEE 754
% \end{keywords}

% \begin{AMS}
% 65F99
% \end{AMS}

\section{Introduction}

Real Givens rotations are one of the most important tools in numerical linear algebra where they are used for solving a wide variety of problems from computing QR factorizations and solving least-squares problems, to finding eigenvalues or computing the singular value decomposition. As a result, many applications from scientific computing to artificial intelligence depend critically on them. Unfortunately, the most common algorithms for constructing them require a square-root operation and many advanced processors do not have a hardware square-root. Consequently there has been substantial interest in square-root free algorithms. The first, designed specifically for the solution of least-squares problems, was proposed by \cite{GENTLEMAN1973} but others followed \cite{Barlow1987,Hammarling1974ANO,Ling1989EfficientLL}. A unified framework for these approaches is presented in \cite{Hsieh1993}. Unfortunately, none of these approaches is directly applicable to the full array of applications that rely on Givens rotations since they do not construct the actual rotation hence our interest in addressing this issue.

It is noteworthy that most of the existing square-root free algorithms for Givens rotations were developed prior to 1990 when the fused multiply-add operation became available in hardware \cite{IBMPOWER1} for the first time on the IBM POWER1 processor. However, the fused multiply-add is commonly available today and we propose an algorithm that can leverage it to directly construct accurate real Givens rotations without the use of square-roots at a speed that is competitive with traditional square-root based algorithms.

\section{Mathematical Preliminaries}

Given a pair of real numbers $f,g \in \Re$, not both zero, we will say that the real unitary matrix\footnote{Real unitary matrices are also commonly called {\em orthogonal matrices}.}
\begin{equation}
   G(f,g) = \frac{1}{\sqrt{f^2+g^2}}\begin{bmatrix}
       f & g \\
       -g & f
   \end{bmatrix}.
\end{equation}
is the Givens rotation associated with $(f,g)$. Note that it is continuous at every point $(f,g)$ in the Cartesian plane except $(0,0)$ where it is undefined. It is common to let
\begin{align}
   r & = \sqrt{f^2+g^2} \label{norm}\\
   c & = \frac{f}{r} \label{cos}\\
   s & = \frac{g}{r}.\label{sin}
\end{align}
and write
\begin{equation}
G(f,g) = \begin{bmatrix}
    c & s \\
    -s & c
\end{bmatrix}
\end{equation}
We note that there is more than one sign convention for Givens rotations. We have chosen this particular form from \cite{Wilkinson,Anderson} wherein the signs of $f$ and $c$ match, and the signs of $g$ and $s$ match. We note that our algorithms can be easily adapted to the sign convention used in \cite{BindelDemmelKahan} where $c \geq 0$ and the sign of $s$ is positive if the signs of $f$ and $g$ match, and negative otherwise.

Givens rotations are a fundamental tool in numerical linear algebra \cite{Wilkinson,GVL,Demmel1997} where they are used to introduce zeros into specific locations in both vectors and matrices. Although there are several approaches to constructing them, our square-root free approach will be based on the following algorithm:\footnote{This is a variation of the well-known and widely used algorithm 4 from \cite{Anderson} (see also algorithm 5.13 in \cite{GVL}).}

\begin{algorithm} c,s = Givens($f,g$)
   \hrule
   \begin{algorithmic}[1]
      \If {$g = 0 \land f = 0$} \Comment{Optional branch for $(0,0)$ input.}
      \State{{\bf return} $\mbox{\tt copysign}(1,f),g$}
      \EndIf
      \If {$|f| \geq |g|$}
         \State{$t = g/f$}
         \State{$c = \mbox{\tt copysign}(1/\mbox{\tt sqrt}(1+t*t),f)$}
         \State{$s = c * t $}
      \Else
         \State{$t = f/g$}
         \State{$s = \mbox{\tt copysign}(1/\mbox{\tt sqrt}(1+t*t),g)$}
         \State{$c = s*t$}
      \EndIf
      \State{{\bf return} $c,s$}
   \end{algorithmic}
\label{Givens}
\end{algorithm}

The {\em optional} branch on lines 1-3 deserves some explanation as it is rather different than those in \cite{Anderson,GVL}. The author believes that this branch is unnecessary and that the algorithm should return $(\Nan,\Nan)$ in the case of a $(0,0)$ input. This reflects the mathematical reality of the singularity at $(0,0)$ and is consistent with the fact that there is no justification for using a rotation to zero out something that is already zero. Removing the optional branch gives exactly that result and yields a small cost savings for the algorithm. The results are otherwise identical and removing unnecessary branching is generally beneficial.

However, we acknowledge that the singularity at $(0,0)$ was historically resolved (see \cite{Wilkinson}) by introducing a branch that assigns the $2\times2$ identity whenever $g = 0$. It was later noted in \cite{Anderson} that this approach fails to preserve continuity at many other points and a simple modification that uses the IEEE {\tt copysign()} utility to address the issue is proposed there. Unfortunately, it does not deal with signed zeros in a consistent manner and is more complex than needed. The approach used above does properly enforce the sign convention even in the presence of signed zeros\footnote{The sign convention from \cite{BindelDemmelKahan} is easily enforced by making simple changes to lines 2,6,and 10. However, there will be a loss of continuity as noted in \cite{BindelDemmelKahan}.} and is only invoked at the single troublesome input. We note that in our tested implementations we use the Julia utility {\tt flipsign()} instead of {\tt copysign()} because it is faster and yields the same result since the argument whose sign is being manipulated is always positive.

We will restrict ourselves to radix-2 IEEE754 compliant arithmetic and assume that round-to-even is used in the event of a tie. We assume that $\ulp(x)$ is a unit in the last place as defined in \cite{Goldberg}. To wit, for any positive real number $x$ in the normal range of the current format we define $\ulp(x) = 2^{e-p+1}$ for any $x \in [2^e,2^{e+1})$ where $p$ is the precision of the floating-point format. For example, in double precision $\ulp(1) = 2^{-52}$. We define the {\em unit roundoff}, which we denote $\bu = \frac{1}{2}\ulp(1)$ so that, for example, in double precision $\bu = 2^{-53}$.

We will find the following lemma useful.
\begin{lemma}
   \label{shiftlemma}
   Let $x$ be a floating-point number and $k$ be a positive integer. If $x \geq k$ and $\ulp(x) \leq 1$ then $x-k$ is also a floating-point number.
\end{lemma}

The proof is simple. Assume the floating-point system has precision $p$. If $\ulp(x) \leq 1$ then $x = d_1d_2...d_s.d_{s+1}...d_p$ where $s \leq p$. Since $x \geq k$ it follows that $d_1d_2...d_s \geq k$ and hence $d_1d_2...d_s - k = \bar{d}_1\bar{d}_2...\bar{d}_s$ can be represented in no more than $s$ digits so that $x-k = \bar{d}_1\bar{d}_2...\bar{d}_s.d_{s+1}...d_p$ which is clearly a precision $p$ floating-point number as well. Note that this lemma is independent of radix.

\section{Computing an Accurate Square-Root Free Givens Rotation}

 We take a very simple {\it approximate and compensate} approach to constructing a square-root free real Givens rotation. We will do so by replacing the reciprocal square-root term in algorithm \ref{Givens} with an appropriate approximation that does not require square roots and use that to generate initial estimates $\bar{c}$ and $\bar{s}$ with algorithm \ref{Givens}. We will then ammend the low accuracy estimates with a compensation scheme. We develop the compensation scheme first as it will inform our choice of approximation strategy.

Although the true quantities $c$ and $s$ satisfy three mathematical conditions -
\begin{eqnarray}
   c^2 + s^2 & = & 1 \label{norm1}\\
   cg - sf & = & 0 \label{orth2}\\
   cf + sg & = & r \label{norm2}
\end{eqnarray}
the computed quantities $\bar{c}$ and $\bar{s}$ may not, no matter how accurately they are computed. A compensation scheme that is based on enforcing equations \ref{norm1} and \ref{orth2} was considered in \cite{borges2024fast}, but it is more costly. We will consider instead a {\it renormalization} scheme that is based only on enforcing equation \ref{norm1}. Mathematically, we will enforce this condition by rescaling $\bar{c}$ and $\bar{s}$ with the factor
$$\frac{1}{\sqrt{\bar{c}^2+\bar{s}^2}}$$
although obviously we must do so without need of a square root. If we let $x = 1 - \bar{c}^2-\bar{s}^2$ then the rescaling factor can be recovered approximately using the Maclaurin series expansion
$$\frac{1}{\sqrt{1-x}} = 1 + \frac{x}{2} + \frac{3x^2}{8} + O(x^3)$$
truncated to second order. This has the added advantage of requiring no divisions when sanely implemented.

Computing $x$ must be done very carefully as it involves {\em significant} cancellation and we will leverage the fused multiply-add operation to do so. We begin by recalling the following algorithm, first proposed by Kahan \cite{KahanQdrtcs}, that computes $ab-cd$ to high-relative precision (error less than $1.5\mbox{\tt ulp}$).
\begin{algorithm} \mbox{\tt abminuscd}(a,b,c,d)
   \hrule
   \begin{algorithmic}
      \State{$tmp = -c * d$}
      \State{{\bf return} $\mbox{\tt fma}(a,b,tmp)-\mbox{\tt fma}(c,d,tmp)$}
   \end{algorithmic}
   \hrule
   \label{Kahanalgo}
\end{algorithm}

We proceed, without loss of generality, by considering the case where $|f| \geq |g|$. In this case we rewrite $1 - \bar{c}^2 - \bar{s}^2$ in the form
\begin{equation}
   (1-|\bar{c}|)(1+|\bar{c}|) - \bar{s}^2. \label{normform}
\end{equation}
If we require that our initial floating-point estimate for $\bar{c}$ be constructed in such a manner that the quantities $(1-|\bar{c}|)$ and $(1+|\bar{c}|)$ are exact floating-point numbers themselves then we could use algorithm \ref{Kahanalgo} to evaluate \ref{normform} to high relative precision. In that case we could use the following algorithm for renormalization (assuming that $|a| \geq |b|$):

\begin{algorithm} renormalize(a,b)
   \hrule
   \begin{algorithmic}[1]
      \State{$x = \mbox{\tt abminuscd}(1-a,1+a,b,b)$}
      \State{$d = x*(1/2+3/8*x)$}
      \State{{\bf return} $(d*a+a),(d*b+b)$}
   \end{algorithmic}
   \label{Renorm}
\end{algorithm}

In order to achieve the highest accuracy it is critical to apply the rescaling as an additive correction (as in line 3 of algorithm \ref{Renorm}) and not in multiplicative form. However, the fused multiply-add is not required for this step.

Continuing under the assumption that $|f| \geq |g|$ we now discuss the specifics of constructing the initial approximation for $\bar{c}$. We will do so indirectly using a simple polynomial or rational approximation $\hat{p}(t)$ to the function
$$
p(t) = 1 + \frac{1}{\sqrt{1+t^2}}
$$
on the interval $-1 \leq t \leq 1$. Note that $p(t)$ is an even function of $t$ and therefore it will suffice to find an approximation $\hat{p}(t)$ that is restricted to the interval $[0,1]$ and evaluate $\hat{p}(|t|)$. We will further require that our approximation satisfy $1 \leq \hat{p}(t) \leq 4$ which is minimally restrictive since $1+ 1/\sqrt{2} \leq p(t) \leq 2$. This restriction guarantees that the quantity $2-\hat{p}(t)$ is an exact floating-point number as a result of Sterbenz lemma (see \cite{FPHB}). Furthermore, $\hat{p}(t) - 1$ is also a floating-point number following lemma \ref{shiftlemma}. If we let $|\bar{c}| = \hat{p}(t) - 1$, then $1+|\bar{c}| = \hat{p}(t)$ and $1-|\bar{c}| = 2-\hat{p}(t)$ which allows us to compute the error in normality to high relative precision which is critical for the renormalization.\footnote{One can accomplish the same thing by approximating $1/\sqrt{1+t^2}$ directly and then setting  the final bit in the significand to zero with a masking operation to generate the initial estimate to $\bar{c}$. This will guarantee that $1+|\bar{c}|$ and  $1-|\bar{c}|$ are both floating-point numbers provided that $1/2 \leq \bar{c} \leq 2$.}

We will use low degree minimax polynomial or rational approximations evaluated in the same precision as the input data (to include the use of appropriately rounded coefficients). Experience shows that the minimax approximations appearing in tables \ref{pcoeffs} and \ref{pcoeffshex} are sufficient to develop our specific algorithms. In particular, we will use the linear approximation as the starting point for Float16 precision computations, the cubic approximation for Float32, and the $2,3$ rational approximation for Float64. These approximations were constructed using the {\tt Chebfun} package in Matlab and all of them take on values well within the interval $[1/2,2]$ for all inputs in $[0,1]$. We note that the maximum absolute error over the interval $[0,1]$ for the linear approximation is roughly $2.3\times10^{-2}$, for the cubic it is $6\times10^{-4}$, and for the rational approximation $6.1\times10^{-7}$.

\begin{table}[p]
   \caption{Coefficients for $\hat{p}(t)$ in decimal form.}
   \label{pcoeffs}
   \begin{tabular}{c l l l l}
      &  &  & \multicolumn{2}{c}{Rational P/Q (Float64)} \\
      & \multicolumn{1}{c}{Linear (Float16)} & \multicolumn{1}{c}{Cubic (Float32)} & \multicolumn{1}{c}{P} & \multicolumn{1}{c}{Q} \\
      \hline
      $a_0$ &  2.023186362360453e+00 &  2.000592060976269e+00 &  4.320202416574703e+01 & 2.160100903948633e+01 \\
      $a_1$ & -2.928932188134525e-01 & -5.865756094031699e-03 &  7.269010078783122e+00 & 3.635125699920760e+00 \\
      $a_2$ &                        & -5.348663062070796e-01 &  2.399473986417590e+01 & 1.738492137035693e+01 \\
      $a_3$ &                        &  2.478388434876587e-01 &                        & 1
   \end{tabular}
\end{table}

\begin{table}[p]
   \caption{Coefficients for $\hat{p}(t)$ rounded to the corresponding precision in hexadecimal form.}
   \label{pcoeffshex}
   \begin{tabular}{c l l l l}
      &  &  & \multicolumn{2}{c}{Rational P/Q (Float64)} \\
      & \multicolumn{1}{c}{Linear (Float16)} & \multicolumn{1}{c}{Cubic (Float32)} & \multicolumn{1}{c}{P} & \multicolumn{1}{c}{Q} \\
      \hline
      $a_0$ &  0x1.03p+1 &  0x1.001366p+1 &  0x1.599dbed88714dp+5 & 0x1.599dbba7931b4p+4 \\
      $a_1$ & -0x1.2cp-2 & -0x1.806b0ep-8 &  0x1.d137760caabecp+2 & 0x1.d14bcc87011f8p+1  \\
      $a_2$ &            & -0x1.11dap-1   &  0x1.7fea74590a9b9p+4 & 0x1.1628a34f936ebp+4 \\
      $a_3$ &            &  0x1.fb92eep-3 &                       & 0x1p+0
   \end{tabular}
\end{table}

Putting it all together yields the square-root free algorithm:

\begin{algorithm} c,s = SqrtFreeGivens(f,g)
   \hrule
   \begin{algorithmic}[1]
      \If {$g = 0$} \Comment{Optional branch for $g = 0$ input.}
         \State{{\bf return} $\mbox{\tt copysign}(1,f),g$}
      \EndIf
      \If  {$|f| \geq |g|$}
         \State{$t = g / f$}
         \State{$\mbox{oneplus} = \hat{p}(|t|)$}
         \State{$c = \mbox{\tt copysign}(\mbox{oneplus}-1,f)$}
         \State{$s = c * t $}
         \State{$\mbox{tmp} = -s*s $}
         \State{$\mbox{err} = \mbox{\tt fma}(\mbox{oneplus},2-\mbox{oneplus},\mbox{tmp}) - \mbox{\tt fma}(s,s,\mbox{tmp})$}
      \Else
         \State{$t = f / g$}
         \State{$\mbox{oneplus} = \hat{p}(|t|)$}
         \State{$s = \mbox{\tt copysign}(\mbox{oneplus}-1,g)$}
         \State{$c = s * t $}
         \State{$\mbox{tmp} = -c*c $}
         \State{$\mbox{err} = \mbox{\tt fma}(\mbox{oneplus},2-\mbox{oneplus},\mbox{tmp}) - \mbox{\tt fma}(c,c,\mbox{tmp})$}
      \EndIf
      \State{$d = \mbox{err}*(1/2 + 3/8*\mbox{err})$}
      \State{{\bf return} $(d*c+c),(d*s+s)$}
   \end{algorithmic}
   \label{RootFreeGivens}
\end{algorithm}

We note that x86 processors include the single-precision {\tt RSQRTSS} instruction that generates a very rapid approximation to the reciprocal square root (roughly the same cost as a multiply \cite{AgnerFog}). In single precision, this approximation has a relative error less than $1.5\times2^{-12}$ (see \cite{Intel}) which is less than that of our cubic minimax approximation and can also be used to generate an appropriate initial approximation very quickly. A version of the square-root free algorithm using this approximation instead of the minimax approximations is included in the package. Experiments show performance and accuracy that is comparable to the version that uses the minimax approximation.

% \section{Error Analysis}

% Without loss of generality we consider the case where $|f| \geq |g|$ and assume that our initial estimate $\bar{c}$ satisfies
% $$\bar{c} = \frac{c}{1+\eps}.$$
% Then the computed value $\bar{s} = \bar{c}*f/g$ satisfies
% $$\bar{s} = \frac{s}{1+\eps}(1+\d1)$$
% for some $|\d1| \leq 2 \bu$. It follows that
% $$
%    1 - \bar{c}^2 - \bar{s}^2 = \frac{(1+\eps)^2-1}{(1+\eps)^2} - s^2\frac{(1+\d1)^2-1}{(1+\eps)^2}
% $$

\section{Testing}

We will use simple simulation to compare the speed and accuracy of algorithm \ref{Givens} versus the square-root-free algorithm \ref{RootFreeGivens}. All testing is done in Julia 1.6 on a Dell Precision desktop computer with an Intel(R) Core(TM) i7-7700K CPU @ 4.20GHz with 16.0 GB of RAM.

We begin by testing performance using the Julia benchmarking utilties with calls of the form:

{\tt @benchmark Givens(f,g) setup=(f=randn(T);g = randn(T))}

\noindent
where {\tt T} is the precision (either Float32 or Float64) being tested.\footnote{We do not do performance testing for Float16 since there is no hardware support for this precision on our processor.} We do this for both algorithms using both precisions. The results of the tests are summarized in table \ref{GivensSpeed} using code that includes the optional branch. We see that the square-root free algorithm runs 20-25\% longer than the square-root based algorithm (running with a hardware square-root).

\begin{table}[p]
   \caption{Median time to compute the Givens rotation with the optional branch.}
   \label{GivensSpeed}
   \begin{tabular}{c c c}
      \toprule
      Precision & Givens & SqrtFreeGivens \\
      \midrule
      Float32   & 3.4ns & 4.2ns \\
      Float64   & 6.5ns & 7.8ns \\
      \toprule
   \end{tabular}
\end{table}

To test the accuracy we use random inputs $f$ and $g$ in the desired precision. We first compute values for $c$ and $s$ by evaluating formulas \ref{norm}, \ref{cos}, and \ref{sin} directly using the greatly extended precision of the BigFloat format in Julia. We round the results of this extended precision calculation to the target precision to get our baseline values. We caution the reader that we are merely comparing to a calculation made with much higher intermediate precision and note that this is not a test for correct rounding. We then compare the outputs of both algorithms against the baseline values to see how well they reproduce the baseline values. We will use $10^{9}$ normally distributed random inputs, that is both $f,g \sim \mathcal{N}(0,1)$. We present the results for all three precisions in table \ref{Floaterror}. We note that the Float16 calculations were done in an emulation mode and not in hardware.

\begin{table}[p]
   \caption{Error rates (\%) with $f,g \sim \mathcal{N}(0,1)$}
   \label{Floaterror}
   \begin{tabular}{l l r r r r}
      \toprule
      & & \multicolumn{2}{c}{Givens} & \multicolumn{2}{c}{SqrtFreeGivens} \\
      Precision & Error & Cosine & Sine & Cosine & Sine \\
      \midrule
      Float 64 & Zero ulp & 57.6 & 57.6 & 82.6 & 82.6 \\
               & One ulp  & 41.3 & 41.3 & 17.4 & 17.4 \\
               & Two ulp  &  1.1 &  1.1 & 0.01 & 0.01 \\
      \toprule
      Float 32 & Zero ulp & 57.7 & 57.7 & 82.6 & 82.6 \\
               & One ulp  & 41.3 & 41.3 & 17.4 & 17.4 \\
               & Two ulp  &  1.0 &  1.0 & 0.01 & 0.01 \\
      \toprule
      Float 16 & Zero ulp & 58.0 & 58.0 & 82.0 & 82.0 \\
               & One ulp  & 41.2 & 41.2 & 17.9 & 17.9 \\
               & Two ulp  &  0.8 &  0.8 & 0.01 & 0.01 \\
      \toprule
   \end{tabular}
\end{table}

\section{Conclusions}

We have developed an accurate square-root free algorithm for constructing real Givens rotations that is based on a well-known and widely used form for constructing them. On a processor that supports a hardware fused multiply-add, the square-root free algorithm runs roughly 20-25\% longer than the square-root based algorithm implemented with a hardware square-root. This makes it a very attractive approach for use on hardware that supports the fused multiply-add but not the square-root. We have also shown experimental evidence that the algorithm is more accurate, on average, than the square-root based approach.

\end{document}